\documentclass[12pt]{article}
\usepackage{geometry}
\usepackage{amssymb,scalerel,amsmath,amsthm}
\usepackage{thmtools, thm-restate}
\usepackage[plainpages = false]{hyperref}
\usepackage[english]{babel}
\usepackage[capitalize]{cleveref}
\usepackage[utf8]{inputenc}
\usepackage{enumerate}
\usepackage{titlesec}
\usepackage[backend = biber, style = numeric, giveninits]{biblatex}

\DeclareMathOperator*{\Times}{\scalerel*{\times}{\textstyle\sum}}

\DeclareMathOperator{\Aut}{Aut}
\DeclareMathOperator{\Diag}{Diag}
\DeclareMathOperator{\End}{End}
\DeclareMathOperator{\Fix}{Fix}

\DeclareMathOperator{\GL}{GL}
\DeclareMathOperator{\Hom}{Hom}
\DeclareMathOperator{\Id}{Id}
\DeclareMathOperator{\im}{Im}

\DeclareMathOperator{\PSL}{PSL}

\DeclareMathOperator{\SL}{SL}

\DeclareMathOperator{\SpecR}{Spec_R}
\DeclareMathOperator{\Spec}{Spec}
\DeclareMathOperator{\Stab}{Stab}

\newcommand*{\eg}{e.g.\ }

\newcommand*{\grpgen}[1]{\left\langle{#1}\right\rangle}

\newcommand*{\ie}{i.e.\ }
\newcommand*{\invb}[1]{\inv{\left(#1\right)}}
\newcommand*{\inv}[1]{{#1}^{-1}}
\newcommand*{\I}{\mathcal{I}}
\newcommand*{\J}{\mathcal{J}}
\newcommand*{\M}{\mathcal{M}}
\newcommand*{\middlebar}{\ \middle | \ }

\newcommand*{\N}{\mathbb{N}}

\newcommand*{\Rconj}[1]{\sim_{#1}}
\newcommand*{\Rinf}{R_{\infty}}

\newcommand*{\Reid}{\mathcal{R}}

\newcommand*{\size}[1]{\left| #1 \right|}

\newcommand*{\Z}{\mathbb{Z}}

\renewcommand{\phi}{\varphi}

\let\originalleft\left
\let\originalright\right
\renewcommand{\left}{\mathopen{}\mathclose\bgroup\originalleft}
\renewcommand{\right}{\aftergroup\egroup\originalright}

\declaretheorem[style=definition, name = Definition, numberwithin=section]{defin}
\declaretheorem[style=definition, name = Example, sibling=defin]{example}
\declaretheorem[name = Theorem, sibling=defin]{theorem}
\declaretheorem[name = Lemma, sibling=defin]{lemma}
\declaretheorem[name = Proposition, sibling=defin]{prop}
\declaretheorem[name = Corollary, sibling=defin]{cor}
\declaretheorem[style=remark, name = Remark, numbered = no]{remark}

\declaretheorem[name = Question, numbered = no]{quest*}

\numberwithin{equation}{section}

\crefname{prop}{Proposition}{Propositions}
\crefname{cor}{Corollary}{Corollaries}

\makeatletter
\renewcommand\maketitle
	{\noindent
		{\Large\bfseries\@title}%
		\bigskip\par\noindent
		{\centering\large\@author}%
		\bigskip
	}

\renewenvironment{abstract}{%
	\par\addvspace{10pt plus2pt minus1pt}\noindent{\bfseries \abstractname\\}\small\leftskip1pc}{\par\vspace{1em}\@endparenv}

\newenvironment{keywords}{%
        \par\addvspace{10pt plus2pt minus1pt}
  \noindent{\bfseries \keywordsname\\}\leftskip1pc\small%
}{%
        \par\@endparenv}

\newenvironment{MSC}{%
        \par\addvspace{10pt plus2pt minus1pt}
  \noindent{\bfseries \mscname\\}\leftskip1pc\small%
}{%
        \par\@endparenv}
\makeatother

\titleformat*{\section}{\large\bfseries}

\title{Twisted Conjugacy in Direct Products of Groups}
\author{Pieter Senden\footnote{KU Leuven Campus Kulak Kortrijk, E.\ Sabbelaan 53, 8500 Kortrijk, Belgium.\\ Contact: pieter.senden@kuleuven.be \, ORCID: 0000-0002-3107-6775}}
\begin{document}
\renewcommand\abstractname{ABSTRACT}
\newcommand{\keywordsname}{KEYWORDS}
\newcommand{\mscname}{2010 Mathematics Subject Classification}
\maketitle
\begin{abstract}
	Given a group \(G\) and an endomorphism \(\phi\) of \(G\), two elements \(x, y \in G\) are said to be \(\phi\)-conjugate if \(x = gy \inv{\phi(g)}\) for some \(g \in G\). The number of equivalence classes for this relation is the Reidemeister number \(R(\phi)\) of \(\phi\). The set \(\{R(\psi) \mid \psi \in \Aut(G)\}\) is called the Reidemeister spectrum of \(G\). We investigate Reidemeister numbers and spectra on direct products of finitely many groups and determine what information can be derived from the individual factors.
\end{abstract}
\begin{MSC}
	20E36, 20E45 (Primary) 20E22 (Secondary)
\end{MSC}
\begin{center}
	This is an Accepted Manuscript of an article published by Taylor \& Francis in Communications in Algebra on 19/07/2021 available at: https://doi.org/10.1080/00927872.2021.1945615.	
\end{center}

\section{Introduction}
Let \(G\) be a group and \(\phi: G \to G\) be an endomorphism. For \(x, y \in G\), we say that \(x\) and \(y\) are \emph{\(\phi\)-conjugate} if there exists a \(g \in G\) such that \(x = g y \inv{\phi(g)}\). In that case, we write \(x \Rconj{\phi} y\) and we denote the \(\phi\)-equivalence class (or Reidemeister class of \(\phi\)) of \(x\) by \([x]_{\phi}\). We also speak of \emph{twisted conjugacy}.

	We define \(\Reid[\phi]\) to be the set of all \(\phi\)-equivalence classes and the \emph{{Reidemeister} number \(R(\phi)\) of \(\phi\)} as the cardinality of \(\Reid[\phi]\). Note that \(R(\phi) \in \N_0 \cup \{\infty\}\). Finally, we define the \emph{Reidemeister spectrum} to be
	\(
		\Spec_R(G) := \{ R(\phi) \mid \phi \in \Aut(G)\}.
	\)
	We say that \(G\) has the \emph{\(\Rinf\)-property} if \(\SpecR(G) = \{\infty\}\). In that case, we also write \(G \in \Rinf\).
	
	The concept of Reidemeister numbers arises from Nielsen fixed-point theory, where the topological analog is used to count and bound the number of fixed-point classes of a continuous self-map, and it is strongly related to the algebraic one introduced above, see \cite{Jiang83}. Besides fixed-point theory, twisted conjugacy also has its applications in representation theory (see \eg \cite{FelshtynLuchnikovTroitsky15,Springer06}), Galois cohomology (see \eg \cite{Serre97}) and isogredience classes (see \eg \cite{FelshtynTroitsky15}).
	
	For many (families of) groups, it has been determined what their Reidemeister spectrum is and/or whether or not they possess the \(\Rinf\)-property, \eg certain subgroups of infinite symmetric groups \cite{Cox19}, extensions of \(\SL(n, \Z)\) and \(\GL(n, \Z)\) by countable abelian groups \cite{MubeenaSankaran14} and Houghton groups \cite{JoLeeLee17}. We refer the reader to \cite{FelshtynNasybullov16} for a more exhaustive list of examples.

	The behavior of the Reidemeister spectrum and the \(\Rinf\)-property under group constructions has been studied as well. One of the major results is due to D.\ Gon\c{c}alves, P.\ Sankaran and P.\ Wong, who proved in \cite[Theorem~1]{GoncalvesSankaranWong20} that under some mild conditions any free product of finitely many (non-trivial) groups has the \(\Rinf\)-property. There are also several results concerning the relation among Reidemeister classes and numbers on group extensions in general, see \eg \cites{Goncalves98,Heath85,Wong01}, and wreath products of abelian groups, see \cite{GoncalvesWong06}.	
	
	In this article, we investigate the following question.
	\begin{quest*}
		Let \(G\) and \(H\) be (non-isomorphic) groups and let \(n \geq 2\) be an integer. What can we say about \(\SpecR(G^{n})\) and \(\SpecR(G \times H)\) in terms of \(\SpecR(G)\) and \(\SpecR(H)\)?
	\end{quest*}
	There has already been done some research into Reidemeister spectra of direct products: S.\ Tertooy determined the Reidemeister spectrum of direct products of free nilpotent groups in his PhD thesis \cite[\S6.5]{Tertooy19} (see also Remark following \cref{cor:SpecRProductCentrelessDirectlyIndecomposableGroups}); K.\ Dekimpe and D.\ Gon\c{c}alves investigated in \cite{DekimpeGoncalves15} the Reidemeister spectrum of infinite abelian groups, some of which are given by an infinite direct product of finite abelian groups. However, these results concern more concrete (families of) groups, whereas we aim to find more general results.
	
	We start by providing a matrix description of the endomorphism monoid of a direct product of groups, which we then use to determine the Reidemeister number of endomorphisms of specific forms. We also derive sufficient conditions to obtain complete information on \(\SpecR(G \times H)\) if we know \(\SpecR(G)\) and \(\SpecR(H)\). We illustrate our results by means of examples.

\section{Matrix Description of Endomorphism Monoid of Direct Product}
Given a group \(G\), the set \(\End(G)\) of all endomorphisms on \(G\) forms a monoid under composition, with the identity map as neutral element. For endomorphisms of direct products of groups, there exists an alternative way to represent this monoid by means of matrices of group homomorphisms, as described by \eg F.\ Johnson in \cite[\S1]{Johnson83}, J.\ Bidwell, M.\ Curran and D.\ McCaughan in \cite[Theorem~1.1]{BidwellCurranMcCaughan06} and by J.\ Bidwell in \cite[Lemma~2.1]{Bidwell08}.

If \(\phi, \psi: G \to H\) are two homomorphisms with commuting images, then it is easily seen that \(\phi + \psi: G \to H: g \mapsto (\phi + \psi)(g) := \phi(g) \psi(g)\) is a homomorphism as well.

Now, let \(G_{1}, \ldots, G_{n}\) be groups. Define
\[
	\M = \left\{ \begin{pmatrix} \phi_{11} & \ldots & \phi_{1n} \\ \vdots & \ddots & \vdots \\ \phi_{n1} & \ldots &\phi_{nn} \end{pmatrix} \middlebar \begin{array}{ccc} \forall 1 \leq i, j \leq n: \phi_{ij} \in \Hom(G_{j}, G_{i}) \\ \forall 1 \leq i, k, l \leq n: k \ne l \implies [\im \phi_{ik}, \im \phi_{il}] = 1\end{array}\right\}
\]
and equip it with matrix multiplication, where the addition of two homomorphisms \(\phi, \psi \in \Hom(G_{j}, G_{i})\) with commuting images is defined as above and the multiplication of \(\phi \in \Hom(G_{j}, G_{i})\) and \(\psi \in \Hom(G_{i}, G_{k})\) is defined as \(\phi \circ \psi\). It is readily verified that this puts a monoid structure on \(\M\), where the diagonal matrix with the respective identity maps on the diagonal is the neutral element.

\begin{lemma}
	For \(G = \Times\limits_{i = 1}^{n}G_{i}\), we have that \(\End (G) \cong \M\) as monoids.
\end{lemma}
F.\ Johnson proved this result for \(G_{1} = \ldots = G_{n}\) and although J. Bidwell, M.\ Curran and D.\ McCaughan state in both aforementioned papers that they only consider finite groups, their proof holds up for infinite groups as well. For the sake of completeness, we give a proof here as well.
\begin{proof}
	For \(1 \leq i \leq n\), denote by \(\pi_{i}: G \to G_{i}\) the canonical projection and by \(e_{i}: G_{i} \to G\) the canonical inclusion. Given \(\phi \in \End (G)\), put \(\phi_{ij} := \pi_{i} \circ \phi \circ e_{j} \in \Hom(G_{j}, G_{i})\).

Fix \(i, k, l \in \{1, \ldots, n\}\) with \(k \ne l\). If \(g_{k} \in G_{k}\) and \(g_{l} \in G_{l}\), then \(e_{k}(g_{k})\) and \(e_{l}(g_{l})\) commute in \(G\). Hence, \(\phi(e_{k}(g_{k}))\) and \(\phi(e_{l}(g_{l}))\) commute as well implying that \(\phi_{ik}(g_{k})\) and \(\phi_{il}(g_{l})\) commute too. Therefore, \([\im \phi_{ik}, \im \phi_{il}] = 1\).

	Since the commuting condition is satisfied, we can define \(F: \End(G) \to \M\) by putting \(F(\phi) := (\phi_{ij})_{ij}\). If \(\phi, \psi \in \End(G)\), then we need to prove for all \(1 \leq i, j \leq n\) that
	\[
		\pi_{i} \circ \phi \circ \psi \circ e_{j} = (\phi \circ \psi)_{ij} = \sum_{k = 1}^{n} \phi_{ik} \psi_{kj}.
	\]
	For \(1 \leq i, j \leq n\) and \(g \in G_{j}\), we see that
	\begin{align*}
		(\pi_{i} \circ \phi \circ \psi \circ e_{j})(g)	&=	(\pi_{i} \circ \phi \circ \psi)(1, \ldots, 1, g, 1, \ldots, 1)	\\
										&=	(\pi_{i} \circ \phi)(\psi_{1j}(g), \ldots, \psi_{nj}(g))	\\
										&=	(\pi_{i} \circ \phi) \left(e_{1}(\psi_{1j}(g))\ldots e_{n}(\psi_{nj}(g)) \right)	\\
										&=	\prod_{k = 1}^{n} (\pi_{i} \circ \phi \circ e_{k})(\psi_{kj}(g))	\\
										&=	\prod_{k = 1}^{n} (\phi_{ik} \circ\psi_{kj})(g)	\\
										&=	\left(\sum_{k = 1}^{n} \phi_{ik} \psi_{kj}\right)(g),
	\end{align*}
	hence the equality holds. Therefore, \(F\) is a monoid homomorphism. It is also clear that \(F\) is injective.
	
	To prove that \(F\) is surjective, let \((\phi_{ij})_{ij} \in \M\) and define
	\[
		\phi: G \to G: (g_{1}, \ldots, g_{n}) \mapsto \left(\prod_{k = 1}^{n}\phi_{1k}(g_{k}), \ldots, \prod_{k = 1}^{n}\phi_{nk}(g_{k})\right).
	\]
	Due to the commuting conditions and the fact that all \(\phi_{ij}\)'s are group homomorphisms, the map \(\phi\) is a well-defined endomorphism of \(G\) and it is clear that \(F(\phi) = (\phi_{ij})_{ij}\).
\end{proof}
We will often identify an endomorphism of \(G\) with its image under \(F\) and write \(\phi = (\phi_{ij})_{ij}\). In a matrix, we denote the identity map with \(1\) and the trivial homomorphism with \(0\).

\begin{lemma}	\label{lem:automorphismOfDirectProductImpliesNormalImages}
	With the notations as above, let \(\phi \in \Aut(G)\). Then
	\begin{enumerate}[(i)]
		\item for all \(1 \leq i \leq n\), \(G_{i}\) is generated by \(\{\im \phi_{ij} \mid 1 \leq j \leq n\}\);
		\item \(\im \phi_{ij}\) is normal in \(G_{i}\) for all \(1 \leq i, j \leq n\).
	\end{enumerate}
\end{lemma}
\begin{proof}
	Suppose \(\phi\) is an automorphism. Then \(\phi\) is surjective. Let \(1 \leq i \leq n\). Then \(G_{i} = \pi_{i}(\phi(G))\). Since \(G\) is generated by \(\{e_{j}(G_{j}) \mid 1 \leq j \leq n\}\), we see that \(G_{i}\) is generated by \(\{(\pi_{i} \circ \phi \circ e_{j})(G_{j}) \mid 1 \leq j \leq n\} = \{\im \phi_{ij} \mid 1 \leq j \leq n\}\). This proves the first item.
	
	For the second, fix \(1 \leq i, j \leq n\). If we pick \(g \in \im \phi_{ij}\) and \(h \in G_{i}\) arbitrary, we can write \(h = xy\) for some \(x \in \im \phi_{ij}\) and \(y \in \grpgen{\{\im \phi_{ik} \mid k \ne j\}}\), as \(G_{i}\) is generated by \(\im \phi_{i1}, \ldots, \im \phi_{in}\) and the images of \(\phi_{ij}\) and \(\phi_{ik}\) commute if \(j \ne k\). Then
	\[
		[g, h] = \inv{g} \invb{xy} g xy = \inv{g} \inv{y} \inv{x} g xy = \inv{g} \inv{x} g x = [g, x] \in \im \phi_{ij}.
	\]
	Consequently, \([\im \phi_{ij}, G_{i}] \leq \im \phi_{ij}\) implying that \(\im \phi_{ij}\) is normal in \(G_{i}\).
\end{proof}
\begin{lemma}	\label{lem:upperTriangularAutomorphisms}
	With the notations as above, suppose that all automorphisms of \(G\) have a matrix representation that is upper triangular, or all of them are lower triangular. Let \(\phi \in \Aut(G)\). Then \(\phi_{ii} \in \Aut(G_{i})\) for each \(i \in \{1, \ldots, n\}\) and \(\phi_{ij} \in \Hom(G_{j}, Z(G_{i}))\) for all \(1 \leq i \ne j \leq n\).
\end{lemma}
\begin{proof}
	We prove the result for upper triangular matrices, the proof for lower triangular is similar.
	
	Let \(\phi \in \Aut(G)\). Let \(\inv{\phi} = (\psi_{ij})_{ij}\). Fix \(i \in \{1, \ldots, n\}\). Since both \(\phi\) and \(\inv{\phi}\) are upper triangular, we find that
	\[
		\Id_{G_{i}} = (\phi \circ \inv{\phi})_{ii} = \phi_{ii} \circ \psi_{ii}
	\]
	and
	\[
		\Id_{G_{i}} = (\inv{\phi} \circ \phi)_{ii} = \psi_{ii} \circ \phi_{ii},
	\]
	showing that \(\phi_{ii}\) must be an automorphism of \(G_{i}\).
	
	Now, let \(1 \leq i, j \leq n\) be indices with \(i \ne j\). Since \(\im \phi_{ij}\) and \(\im \phi_{ii} = G_{i}\) commute, we conclude that \(\im \phi_{ij} \in Z(G_{i})\).
\end{proof}

Using this alternative description of the endomorphism monoid, we can deduce some general results regarding Reidemeister numbers of specific endomorphisms on direct products. 

We define the diagonal endomorphisms to be all endomorphisms of \(\End (G)\) of the form
\[
	\Diag(\phi_{1}, \ldots, \phi_{n}): G \to G: (g_{1}, \ldots, g_{n}) \mapsto (\phi_{1}(g_{1}), \ldots, \phi_{n}(g_{n})),
\]
where each \(\phi_{i} \in \End (G_{i})\). We denote the submonoid of all diagonal endomorphism with \(\Diag(G)\). Note that it is isomorphic with \(\End(G_{1}) \times \ldots \times \End(G_{n})\).

The following is then quite straightforward.
\begin{prop}	\label{prop:ReidemeisterNumberDirectProductAutomorphismGroups}
	Let $G_{1}, \ldots, G_{n}$ be groups and put \(G = \Times\limits_{i = 1}^{n} G_{i}\). Let \(\phi\) be an element of \(\Diag(G)\) and write \(\phi = \Diag(\phi_{1}, \ldots, \phi_{n})\). Then $R(\phi) = \prod_{i = 1}^{n} R(\phi_{i})$.
\end{prop}
\begin{proof}
	It is clear that, for $(g_{1}, \ldots, g_{n}), (h_{1}, \ldots, h_{n}) \in G_{1} \times \ldots \times G_{n}$, we have
	\[
		(g_{1}, \ldots, g_{n}) \Rconj{\phi} (h_{1}, \ldots, h_{n}) \iff \forall i \in \{1, \ldots, n\}: g_{i} \Rconj{\phi_{i}} h_{i}.
	\]
	Thus, the map
	\[
		\Reid[\phi] \to \Reid[\phi_{1}] \times \ldots \times \Reid[\phi_{n}]: [(g_{1}, \ldots, g_{n})]_{\phi} \mapsto ([g_{1}]_{\phi_{1}}, \ldots, [g_{n}]_{\phi_{n}})
	\]
	is a well-defined bijection, implying that
	\[
		R(\phi) = \prod_{i = 1}^{n} R(\phi_{i}). \qedhere
	\]
\end{proof}
\begin{defin}
	For \(a \in \N_{0} \cup \{\infty\}\), we define the product \(a \cdot \infty\) to be equal to \(\infty\). Let \(A_{1}, \ldots, A_{n}\) be subsets of \(\N_{0} \cup \{\infty\}\). We then define
	\[
		A_{1} \cdot \ldots \cdot A_{n} := \prod_{i = 1}^{n} A_{i} := \{a_{1} \ldots a_{n} \mid \forall i \in \{1, \ldots, n\}: a_{i} \in A_{i}\}.
	\]
	If \(A_{1} = \ldots = A_{n} =: A\), we also write \(A^{(n)}\) for the \(n\)-fold product of \(A\) with itself.
\end{defin}

\begin{cor}	\label{cor:SpecRDirectProductCharacteristicSubgroups}
	Let \(G_{1}, \ldots, G_{n}\) be groups and put \(G = \Times\limits_{i = 1}^{n} G_{i}\). Then
	\[
		\prod_{i = 1}^{n} \SpecR(G_{i}) \subseteq \SpecR(G).
	\]
	Equality holds if \(\Aut(G) = \Times\limits_{i = 1}^{n} \Aut(G_{i})\).
\end{cor}
We now specify to the case of an \(n\)-fold direct product of a group with itself, \ie \(G^{n}\) for some group \(G\) and integer \(n \geq 1\). First of all, the symmetric group \(S_{n}\) embeds in \(\End(G^{n})\) in the following way:
\[
	S_{n} \to \End G^{n}: \sigma \mapsto (P_{\inv{\sigma}}: G \to G: (g_{1}, \ldots, g_{n}) \mapsto (g_{\inv{\sigma}(1)}, \ldots, g_{\inv{\sigma}(n)})),
\]
where the matrix representation of \(P_{\inv{\sigma}}\) is given by
\[
	P_{\inv{\sigma}} = \begin{pmatrix}
		e_{\inv{\sigma}(1)}	\\
		\vdots		\\
		e_{\inv{\sigma}(n)}
	\end{pmatrix}.
\]
Here, \(e_{i}\) is a row with a \(1\) on the \(i\)-th spot and zeroes elsewhere.

To prove that is indeed a monoid morphism, let \(\sigma, \tau \in S_{n}\) and \((g_{1}, \ldots, g_{n}) \in G^{n}\) be arbitrary. Put \(h_{i} := g_{\inv{\sigma}(i)}\), then
\[
	P_{\inv{\tau}}(h_{1}, \ldots, h_{n}) = (h_{\inv{\tau}(1)}, \ldots, h_{\inv{\tau}(n)}) = (g_{\inv{\sigma}(\inv{\tau}(1))}, \ldots, g_{\inv{\sigma}(\inv{\tau}(n))})
\]
and
\[
	P_{\inv{(\tau \sigma)}}(g_{1}, \ldots, g_{n}) = (g_{\inv{\sigma}(\inv{\tau}(1))}, \ldots, g_{\inv{\sigma}(\inv{\tau}(n))}),
\]
therefore \(P_{\inv{(\tau\sigma)}} = P_{\inv{\tau}} P_{\inv{\sigma}}\).

Now, we define \(\End_{w}(G^{n})\) to be the submonoid of \(\End(G^{n})\) generated by \(S_{n}\) and \(\Diag(G^{n})\).

\begin{lemma}
	Let \(G\) be a group and \(n \geq 1\) an integer. Then each endomorphism \(\phi\) in \(\End_{w}(G^{n})\) can be written as
	\[
		\phi = \Diag(\phi_{1}, \ldots, \phi_{n}) P_{\inv{\sigma}}
	\]
	for some \(\phi_{i} \in \End(G)\) and \(\sigma \in S_{n}\). Moreover, \(\phi \in \Aut(G^{n})\) if and only if each \(\phi_{i} \in \Aut(G)\).
\end{lemma}
\begin{proof}
	The existence of such a decomposition relies on the following equality, which we claim holds for all \(\sigma \in S_{n}\) and \(\phi_{i} \in \End(G)\):
	\begin{equation}	\label{eq:rewritingRulePermutationDiagonalMatrices}
		P_{\inv{\sigma}} \Diag(\phi_{1}, \ldots, \phi_{n}) = \Diag(\phi_{\inv{\sigma}(1)}, \ldots, \phi_{\inv{\sigma}(n)}) P_{\inv{\sigma}}.
	\end{equation}
	Indeed, evaluating the left-hand side in \((g_{1}, \ldots, g_{n})\) yields
	\[
		P_{\inv{\sigma}} (\phi_{1}(g_{1}), \ldots, \phi_{n}(g_{n})) = (\phi_{\inv{\sigma}(1)}(g_{\inv{\sigma}(1)}), \ldots, \phi_{\inv{\sigma}(n)}(g_{\inv{\sigma}(n)}))
	\]
	whereas the right-hand side yields
	\[
		\Diag(\phi_{\inv{\sigma}(1)}, \ldots, \phi_{\inv{\sigma}(n)}) (g_{\inv{\sigma}(1)}, \ldots, g_{\inv{\sigma}(n)}),
	\]
	hence we see they are equal.
	
	Thus, given an element in \(\End_{w}(G^{n})\), we can apply the equality above several times in order to gather all diagonal endomorphisms and all elements of the form \(P_{\inv{\sigma}}\) together, yielding the desired representation.
	
	The claim regarding the automorphism is immediate, as each \(P_{\inv{\sigma}}\) is an automorphism.
\end{proof}
The representation is not necessarily unique, as the trivial endomorphism equals \(\Diag(0, \ldots, 0) P_{\inv{\sigma}}\) for all \(\sigma \in S_{n}\). If we restrict ourselves to automorphisms, however, this yields the injective group homomorphism
	\[
		\Psi: \Aut(G) \wr S_{n} \to \Aut(G^{n}): (\phi, \sigma) = (\phi_{1}, \ldots, \phi_{n}, \sigma) \mapsto \Diag(\phi_{1}, \ldots, \phi_{n}) P_{\inv{\sigma}}.
	\]
Here, \(\Aut(G) \wr S_{n}\) is the wreath product, \ie the semidirect product \(\Aut(G)^{n} \rtimes S_{n}\), where the action is given by
\[
	\sigma \cdot (\phi_{1}, \ldots, \phi_{n}) = (\phi_{\inv{\sigma}(1)}, \ldots, \phi_{\inv{\sigma}(n)}).
\]
To see that is indeed homomorphism, note that, for \((\phi, \sigma), (\psi, \tau) \in \Aut(G) \wr S_{n}\),
\[
	\Psi((\phi, \sigma) (\psi, \tau)) = \Psi(\phi \circ (\sigma \cdot \psi), \sigma \tau) = \Diag(\phi_{1}\psi_{\inv{\sigma}(1)}, \ldots, \phi_{n}\psi_{\inv{\sigma}(n)}) P_{\invb{\sigma \tau}}
\]
and
\begin{align*}
	\Psi(\phi, \sigma)\Psi(\psi, \tau)	&=	\Diag(\phi_{1}, \ldots, \phi_{n}) P_{\inv{\sigma}} \Diag(\psi_{1}, \ldots, \psi_{n}) P_{\inv{\tau}}	\\
							&=	\Diag(\phi_{1}, \ldots, \phi_{n}) \Diag(\psi_{\inv{\sigma}(1)}, \ldots, \psi_{\inv{\sigma}(n)}) P_{\inv{\sigma}} P_{\inv{\tau}}	\\
							&=	\Diag(\phi_{1}\psi_{\inv{\sigma}(1)}, \ldots, \phi_{n}\psi_{\inv{\sigma}(n)}) P_{\invb{\sigma \tau}}
\end{align*}
where we used \eqref{eq:rewritingRulePermutationDiagonalMatrices}.
We will identify \(\Aut(G) \wr S_{n}\) with its image under \(\Psi\) and thus regard it as a subgroup of \(\Aut(G^{n})\).

We now determine the Reidemeister number of an arbitrary element of \(\End_{w}(G^{n})\), which also generalizes \cite[Proposition~5.1.2]{DekimpeSenden21}.
\begin{prop}[{see \eg \cite[Proposition~1.1.3]{DekimpeSenden21}}]	\label{prop:conjugateEndomorphisms}
    	Let \(G\) be a group, \(\phi \in \End(G)\) and \(\psi \in \Aut(G)\). Then \(R(\phi) = R(\phi^\psi)\).
\end{prop}

\begin{prop}	\label{prop:ReidemeisterNumberPermutedDiagonalEndomorphism}
	Let \(G\) be a group and let \(\phi := \Diag(\phi_{1}, \ldots, \phi_{n}) P_{\sigma} \in \End_{w}(G^{n})\). Let
	\[
		\sigma = (c_{1} \ldots c_{n_{1}})(c_{n_{1} + 1} \ldots c_{n_{2}})\ldots (c_{n_{k - 1} + 1} \ldots c_{n_{k}})
	\]
	be the disjoint cycle decomposition of \(\sigma\), where \(n_{0} := 0 < n_{1} < n_{2} < \ldots < n_{k} = n\). Put \(\tilde{\phi}_{j} := \phi_{c_{n_{j - 1} + 1}} \circ \ldots \circ \phi_{c_{n_{j}}}\). Then
	\[
		R(\phi) = \prod_{j = 1}^{k} R(\tilde{\phi}_{j}).
	\]

\end{prop}
\begin{proof}
	Let \(\tau \in S_{n}\) be the permutation given by \(\tau(i) = c_{i}\) for all \(i \in \{1, \ldots, n\}\). Then
	\[
		\sigma^{\tau} = \inv{\tau} \sigma \tau = (1 \ldots n_{1})(n_{1} + 1 \ \ldots n_{2}) \ldots (n_{k - 1} + 1 \ldots n_{k}).
	\]
	By \cref{prop:conjugateEndomorphisms}, conjugate endomorphisms have the same Reidemeister number. Since, by \eqref{eq:rewritingRulePermutationDiagonalMatrices}, 
	\[
		P_{\tau}\Diag(\phi_{1}, \ldots, \phi_{n}) P_{\sigma} P_{\inv{\tau}} = \Diag(\phi_{\tau(1)}, \ldots, \phi_{\tau(n)})P_{\sigma^{\tau}},
	\]
	it is sufficient to determine the Reidemeister number of the latter endomorphism. Note that \(P_{\sigma^{\tau}}\) is a block matrix consisting of \(k\) square blocks of the form
	\[
		\begin{pmatrix}
			0		&	1		&	0		&	0		&	\ldots	&	0	&	0		\\
			0		&	0		&	1		&	0		&	\ldots	&	0	&	0		\\
			0		&	0		&	0		&	1		&	\ldots	&	0	&	0		\\
			\vdots	&	\vdots	&	\vdots	&	\vdots	&	\ddots	&	\vdots&	\vdots	\\
			0		&	0		&	0		&	0		&	\ldots	&	1	&	0		\\
			0		&	0		&	0		&	0		&	\ldots	&	0	&	1		\\
			1		&	0		&	0		&	0		&	\ldots	&	0	&	0	
		\end{pmatrix}
	\]
	on the diagonal and zero matrices elsewhere. This implies that
	\[
		\Diag(\phi_{\tau(1)}, \ldots, \phi_{\tau(n)})P_{\sigma^{\tau}} \in \Times_{i = 1}^{k} \End(G^{n_{i} - n_{i - 1}}).
	\]
	Therefore, it is sufficient to determine Reidemeister numbers of endomorphisms of the form \(\psi := \Diag(\psi_{1}, \ldots, \psi_{n})P_{\alpha}\), where \(\alpha = (1 \ 2 \ldots \ n)\). Indeed, if we know that \(R(\psi) = R(\psi_{1} \circ \ldots \circ \psi_{n})\), then
	\begin{align*}
		R(\phi)	&=	\prod_{i = 1}^{k} R(\phi_{\tau(n_{i - 1} + 1)} \circ \ldots \circ \phi_{\tau(n_{i})}) = \prod_{i = 1}^{k} R(\phi_{c_{n_{i - 1} + 1}} \circ \ldots \circ \phi_{c_{n_{i}}}) = \prod_{j = 1}^{k} R(\tilde{\phi}_{j}),
	\end{align*}
	by construction of \(\tau\).
	
	So, let \((g_{1}, \ldots, g_{n}) \in G^{n}\). We claim there exists a \(g \in G\) such that
	\[
		(g_{1}, \ldots, g_{n}) \Rconj{\psi} (g, 1, \ldots, 1).
	\]
	Note that, for \((x_{1}, \ldots, x_{n}) \in G^{n}\),
	\[
		(x_{1}, \ldots, x_{n}) (g_{1}, \ldots, g_{n}) \inv{\psi(x_{1}, \ldots, x_{n})} = (x_{1}g_{1} \inv{\psi_{1}(x_{2})}, \ldots, x_{n}g_{n} \inv{\psi_{n}(x_{1})}).
	\]
	Put \(x_{1} = 1, x_{n} = \inv{g}_{n}\) and \(x_{i} = \psi_{i}(x_{i + 1})\inv{g}_{i}\) for \(i \in \{2, \ldots, n - 1\}\), starting with \(i = n - 1\). Finally, put \(g = x_{1} g_{1} \inv{\psi_{1}(x_{2})}\). Then
	\[
		\begin{cases}
			g = x_{1} g_{1} \inv{\psi_{1}(x_{2})} \\
			1 = x_{i} g_{i} \inv{\psi_{i}(x_{i + 1})}	& \mbox{for \(i \geq 2\)},
		\end{cases}
	\]
	where \(x_{n + 1} = x_{1}\). This implies that \((g_{1}, \ldots, g_{n}) \Rconj{\psi} (g, 1, \ldots, 1)\).
	
	Next, put \(\tilde{\psi} = \psi_{1} \circ \ldots \circ \psi_{n}\) and suppose \((g, 1, \ldots, 1) \Rconj{\psi} (h, 1, \ldots, 1)\) for some \(g, h \in G\). Then there exist \(x_{1}, \ldots, x_{n} \in G\) such that
	\[
		\begin{cases}
			g = x_{1} h\inv{\psi_{1}(x_{2})} \\
			1 = x_{i} \inv{\psi_{i}(x_{i + 1})}	& \mbox{for \(i \geq 2\)},
		\end{cases}
	\]
	where again \(x_{n + 1} = x_{1}\). Consequently, \(x_{i} = \psi_{i}(x_{i + 1})\) for \(i \geq 2\). This implies that
	\[
		g = x_{1}h \inv{\psi_{1}(x_{2})} = x_{1}h \inv{\psi_{1}(\psi_{2}(x_{3}))} = \ldots =  x_{1}h \inv{\tilde{\psi}(x_{1})},
	\]
	\ie \(g \Rconj{\tilde{\psi}} h\). Conversely, if \(g \Rconj{\tilde{\psi}} h\), then there exists an \(x \in G\) such that \(g = x h \inv{\tilde{\psi}(x)}\). Put \(x_{1} = x\) and \(x_{i} = \psi_{i}(x_{i + 1})\) for \(i \geq 2\), starting with \(i = n\) and where, again, \(x_{n + 1} = x_{1}\). Then
	\[
		\begin{cases}
			g = x_{1} h\inv{\psi_{1}(x_{2})} \\
			1 = x_{i} \inv{\psi_{i}(x_{i + 1})}	& \mbox{for \(i \geq 2\)},
		\end{cases}
	\]
	hence \((g, 1, \ldots, 1) \Rconj{\psi} (h, 1, \ldots, 1)\).
	
	Combining all results, we find that there is a bijection between the Reidemeister classes of \(\psi\) and those of \(\tilde{\psi}\). Consequently, \(R(\psi) = R(\tilde{\psi})\).
\end{proof}
\begin{cor}	\label{cor:ReidemeisterNumbersOfAut(G)wrSn}
	Let \(G\) be a group and \(n \geq 1\) a natural number. Then
	\[
		\bigcup_{i = 1}^{n} \SpecR(G)^{(i)} \subseteq \SpecR(G^{n}).
	\]
	Equality holds if \(\Aut(G^{n}) = \Aut(G) \wr S_{n}\).
\end{cor}
\begin{proof}
	Let \(1 \leq k \leq n\). Let \(\phi_{1}, \ldots, \phi_{k}\) be automorphisms of \(G\). We prove that \(R(\phi_{1}) \ldots R(\phi_{k}) \in \SpecR(G^{n})\). Consider the automorphism
	\[
		\phi := \Diag(\phi_{1}, \ldots, \phi_{k}, \Id_{G}, \ldots, \Id_{G})P_{\sigma}
	\]
	where \(\sigma = (1)(2)(3) \ldots (k - 1)(k \ k + 1 \ \ldots \ n - 1 \ n)\). By the previous proposition \(R(\phi) = R(\phi_{1}) \ldots R(\phi_{k})\).
	
	Note that this combined with the previous proposition also proves that the left-hand side equals \(\{R(\phi) \mid \phi \in \Aut(G) \wr S_{n} \}\), from which the additional claim follows immediately.
\end{proof}
\begin{cor}	\label{cor:RinfWhenAut(Gn)=Aut(G)wrSn}
	Let \(G\) be a group and \(n \geq 1\) an integer. Suppose that \(G \in \Rinf\) and that \(\Aut(G^{n}) = \Aut(G) \wr S_{n}\). Then \(G^{n} \in \Rinf\).
\end{cor}
\begin{proof}
	Since \(\Aut(G) \wr S_{n} = \Aut(G^{n})\), the previous corollary shows that
	\[
		\SpecR(G^{n}) = \bigcup_{i = 1}^{n} \SpecR(G)^{(i)} = \bigcup_{i = 1}^{n} \{\infty\}^{(i)} = \{\infty\}. \qedhere
	\]
\end{proof}
\begin{example}
	Although \cref{cor:ReidemeisterNumbersOfAut(G)wrSn} yields some information regarding \(\SpecR(G^{n})\), this information can be limited. Consider the case \(G = \Z\). It is well-known (see \eg \cite{Romankov11}) that, for \(r \geq 1\),
	\[
		\SpecR(\Z^{r}) = \begin{cases}
			\{2, \infty\}		&	\mbox{if } r = 1	\\
			\N_{0} \cup \{\infty\}	&	\mbox{otherwise.}
		\end{cases}
	\]
	However, using only the result for \(r = 1\), \cref{cor:ReidemeisterNumbersOfAut(G)wrSn} merely yields \(\{2^{i} \mid 1 \leq i \leq n\} \cup \{\infty\} \subseteq \SpecR(\Z^{n})\). 
	Nonetheless, we will later provide conditions under which \cref{cor:ReidemeisterNumbersOfAut(G)wrSn,cor:RinfWhenAut(Gn)=Aut(G)wrSn} yield full information on \(\SpecR(G^{n})\).
\end{example}
We can also generalize \cref{prop:conjugateEndomorphisms}.
\begin{cor}	\label{cor:CyclicPermutationsPreserveReidemeisterNumber}
	Let \(G\) be a group, \(n \geq 1\) an integer and \(\phi_{1}, \ldots, \phi_{n} \in \End(G)\). Then
	\[
		R(\phi_{1} \circ \phi_{2} \circ \ldots \circ \phi_{n}) = R(\phi_{2} \circ \phi_{3} \ldots \circ \phi_{n} \circ \phi_{1})
	\]
\end{cor}\begin{proof}
	Consider the endomorphism \(\phi := \Diag(\phi_{1}, \ldots, \phi_{n}) P_{(1 \ 2 \ \ldots \ n)}\) of \(G^{n}\). Since \((1 \ 2 \ \ldots \ n)\) and \((2 \ \ldots \ n \ 1)\) are both cycle representations of the same permutation, \cref{prop:ReidemeisterNumberPermutedDiagonalEndomorphism} yields that
	\[
		R(\phi_{1} \circ \phi_{2} \circ \ldots \circ \phi_{n}) = R(\phi) = R(\phi_{2} \circ \phi_{3} \circ \ldots \circ \phi_{n} \circ \phi_{1}). \qedhere
	\]
\end{proof}

\section{Direct Products of Two Groups}
We now restrict to the case of direct products of two groups. Let \(H\) and \(K\) be two groups and put \(G := H \times K\). Instead of using indices to represent endomorphisms of \(G\) as matrices, we use the notation
\[
	\M = \left\{ \begin{pmatrix} \alpha & \beta \\ \gamma & \delta \end{pmatrix} \middlebar \begin{array}{ccc}\alpha \in \End(H), & \beta \in \Hom(K, H), &[\im \alpha, \im \beta] = 1\\ \gamma \in \Hom(H, K), &\delta \in \End(K), &[\im \gamma, \im \delta] = 1\end{array}\right\}.
\]

Let \(\phi \in \End(G)\) be an endomorphism leaving \(H\) invariant. Then, in the notation above, \(\gamma\) is the trivial homomorphism, sending everything to \(1\). Fix representatives \(\{h_{i}\}_{i \in \mathcal{I}}\) of \(\Reid[\alpha]\) and \(\{k_{j}\}_{j \in \mathcal{J}}\) of \(\Reid[\delta]\).

\begin{lemma}	\label{lem:representativesReidphi}
	If \(\phi\) is an endomorphism of \(G = H \times K\) leaving \(H\) invariant, then with the notation as above, we have \(\Reid[\phi] =  \{[(h_{i},k_{j})]_{\phi} \mid i \in \mathcal{I}, j \in \mathcal{J}\}\). In particular, \(R(\phi) \leq R(\alpha)R(\delta)\).
\end{lemma}
\begin{proof}
	Let \((h, k) \in G\). We have to prove that there are indices \(i\) and \(j\) such that \((h, k) \Rconj{\phi} (h_{i}, k_{j})\). First, let \(j \in \mathcal{J}\) be such that \(k \Rconj{\delta} k_{j}\). Write \(k_{j} = y k \inv{\delta(y)}\) for some \(y \in K\), then
	\[
		(1, y) (h, k) \inv{\phi(1, y)} = (h \inv{\beta(y)}, k_{j}).
	\]
	Put \(h' = h \inv{\beta(y)}\), then \((h, k) \Rconj{\phi} (h', k_{j})\). Next, let \(i \in \mathcal{I}\) be such that \(h' \Rconj{\alpha} h_{i}\) and write \(h_{i} = x h' \inv{\alpha(x)}\) for some \(x \in H\). Then
	\[
		(x, 1)(h', k_{j})\inv{\phi(x, 1)} = (h_{i}, k_{j}),
	\]
	finishing the proof.
	
	The inequality \(R(\phi) \leq R(\alpha) R(\delta)\) then follows immediately.
\end{proof}
The \(\phi\)-conjugacy relation can also be seen as a group action of \(G\) on itself, namely \(g \cdot h := g h \inv{\phi(g)}\). This allows us to speak of orbits and stabilizers.
\begin{defin}
	Let \(A\) be a group, \(\phi \in \End(A)\) and \(a \in A\). The \emph{\(\phi\)-stabilizer of \(a\)} is the subgroup
	\[
		\Stab_{\phi}(a) = \{b \in A \mid a = b a \inv{\phi(b)}\}
	\]
	of \(A\). We also call subgroups of this form \emph{twisted stabilizers}.
\end{defin}
We continue with the same notation as before. Fix \(j \in \mathcal{J}\). The subgroup \(\Stab_{\delta}(k_{j})\) acts on the right on \(\Reid[\alpha]\) in the following way:
\[
	\rho_{j}: \Reid[\alpha] \times \Stab_{\delta}(k_{j}) \to \Reid[\alpha]: ([h]_{\alpha}, y) \mapsto [h\beta(y)]_{\alpha}.
\]
Indeed, if \(h' = xh \inv{\alpha(x)}\), then
\[
	h' \beta(y) = xh \inv{\alpha(x)} \beta(y) = x h \beta(y) \inv{\alpha(x)},
\]
since \([\im \alpha, \im \beta] = 1\). Therefore, the action is independent of the representative. Moreover, if \(y, y' \in \Stab_{\delta}(k_{j})\), then
\[
	\rho_{i}([h \beta(y')]_{\alpha}, y) = [h \beta(y')\beta(y)]_{\alpha} = [h \beta(y'y)]_{\alpha} = \rho_{i}([h]_{\alpha}, y'y).
\]

\begin{theorem}	\label{theo:sumFormulaReidemeisterNumberDirectProductInvariantFactor}
	With the notations as above, denote by \(r_{j}\) the number of orbits of \(\rho_{j}\). Then
	\[
		R(\phi) = \sum_{j \in \mathcal{J}} r_{j}.
	\]
	Here, an infinite sum or a sum with one of its terms equal to \(\infty\) is to be interpreted as \(\infty\).
\end{theorem}
\begin{proof}
	By \cref{lem:representativesReidphi}, we have to decide when \((h_{i_{1}}, k_{j_{1}})\) and \((h_{i_{2}}, k_{j_{2}})\) are \(\phi\)-conjugate. If so, then, since \(H\) is \(\phi\)-invariant, projecting onto \(K\) yields \(k_{j_{1}} \Rconj{\delta} k_{j_{2}}\), hence \(k_{j_{1}} = k_{j_{2}} =: k\). We claim that
	\[
		(h_{i_{1}}, k) \Rconj{\phi} (h_{i_{2}}, k) \iff \exists (x, y) \in H \times \Stab_{\delta}(k): h_{i_{1}} = x h_{i_{2}} \inv{\alpha(x)} \inv{\beta(y)}.
	\]
	Indeed, if \((h_{i_{1}}, k) = (x, y)(h_{i_{2}}, k) \inv{\phi(x, y)}\) for some \((x, y) \in G\), then
	\begin{align*}
		h_{i_{1}} &= x h_{i_{2}}\inv{\alpha(x)} \inv{\beta(y)}	\\
		k &= y k \inv{\delta(y)},
	\end{align*}
	showing that \(y \in \Stab_{\delta}(k)\). The converse implication is clear. As \([\im \alpha, \im \beta] = 1\), we can rewrite this as
	\[
		(h_{i_{1}}, k) \Rconj{\phi} (h_{i_{2}}, k) \iff \exists y \in \Stab_{\delta}(k): [h_{i_{1}}]_{\alpha} = [h_{i_{2}}\inv{\beta(y)}]_{\alpha}
	\]
	\ie if and only if \(h_{i_{1}}\) and \(h_{i_{2}}\) lie in the same orbit under the action of \(\Stab_{\delta}(k)\). Therefore, the theorem is proven.
\end{proof}
Of course, the analogous result for \(\phi\) leaving \(K\) invariant holds as well.

We will now use this theorem to determine the Reidemeister spectrum of direct products of the form \(G \times F\), where \(F\) is a finite group and \(G\) is a finitely generated torsion-free residually finite group, generalising a result due to A.\ Fel'shtyn (see \cite[Proposition~3]{Felshtyn00}). The key ingredient will be proving that automorphisms of \(G\) with finite Reidemeister number have trivial twisted stabilizers.

The following two results are well-known, see \eg \cite[Corollary~2.5]{FelshtynTroitsky07} and \cite[Lemma~1.1]{GoncalvesWong09}, respectively.
\begin{lemma}	\label{lem:ReidemeisterNumberInvariantInnerAutomorphisms}
	Let \(G\) be a group, \(\phi \in \End(G)\) and \(g \in G\). Denote by \(\tau_{g}\) the inner automorphism corresponding to \(g\), \ie \(\tau_{g}(x) = gx\inv{g}\) for all \(x \in G\). Then \(R(\tau_{g} \circ \phi) = R(\phi)\).
\end{lemma}
\begin{lemma}	\label{lem:ReidemeisterNumbersExactSequence}
	Let \(G\) be a group, \(\phi \in \End(G)\) and \(N\) a \(\phi\)-invariant normal subgroup of \(G\), \ie \(\phi(N) \leq N\). Denote by \(\bar{\phi}\) the induced endomorphism on \(G / N\) and by \(\phi'\) the induced endomorphism on \(N\). Then the following hold:
	\begin{enumerate}[(i)]
		\item \(R(\phi) \geq R(\bar{\phi})\). In particular, if \(R(\bar{\phi}) = \infty\), then also \(R(\phi) = \infty\).	\label{lemitem:ReidemeisterNumberGreaterThanQuotientReidemeisterNumber}
		\item If \(R(\phi') = \infty\) and \(\size{\Fix(\bar{\phi})} < \infty\), then \(R(\phi) = \infty\).	\label{lemitem:InfiniteSubgroupReidemeisterNumberFiniteFixedPointsImpliesInfiniteReidemeisterNumber}
	\end{enumerate}
	In particular, if \(N\) is characteristic and \(G / N\) has the \(\Rinf\)-property, then so does \(G\).
\end{lemma}

We also need the following bound.
\begin{lemma}[{\cite[Lemma~3]{Jabara08}}]	\label{lem:upperboundFixedPointsFiniteGroupReidemeisterNumber}
	Let \(G\) be a finite group and \(\phi \in \Aut(G)\). Put \(R(\phi) = r\). Then \(\size{\Fix(\phi)} \leq 2^{2^{r}}\).
\end{lemma}
The following result can be found implicitly in \cite{Jabara08}. For the reader's convenience, we present it here with a complete proof.

\begin{prop}	\label{prop:infiniteTwistedStabilisersImpliesInfiniteReidemeisterNumberFGRF}
	Let \(G\) be a finitely generated residually finite group. Let \(\phi \in \Aut(G)\). If \(\Stab_{\phi}(g)\) is infinite for some \(g \in G\), then \(R(\phi) = \infty\).
\end{prop}
\begin{proof}
	Suppose that \(\Stab_{\phi}(g)\) is infinite. Note that
	\begin{align*}
		\Stab_{\phi}(g)	&=	\{x \in G \mid x g \inv{\phi(x)} = g\}	\\
					&=	\{x \in G \mid \phi(x) = x^{g}\}	\\
					&=	\Fix(\tau_{g} \circ \phi).
	\end{align*}
	As \(\phi\) and \(\tau_{g} \circ \phi\) have the same Reidemeister number by \cref{lem:ReidemeisterNumberInvariantInnerAutomorphisms}, it is sufficient to prove the result for \(g = 1\).
	
	So, suppose \(\phi\) has infinitely many fixed points. Fix \(n \geq 1\) and let \(x_{1}, \ldots, x_{n}\) be \(n\) of these fixed points. Since \(G\) is finitely generated and residually finite, we can find a characteristic subgroup \(K\) of \(G\) with finite index such that \(\pi: G \to G / K\) is injective on \(\{x_{1}, \ldots, x_{n}\}\). Let \(\bar{\phi}\) be the induced automorphism on \(G / K\). Then \(R(\phi) \geq R(\bar{\phi})\). Moreover, by \cref{lem:upperboundFixedPointsFiniteGroupReidemeisterNumber}, we know that
	\[
		n = \size{\{\bar{x}_{1}, \ldots, \bar{x}_{n}\}} \leq \size{\Fix(\bar{\phi})} \leq 2^{2^{R(\bar{\phi})}},
	\]
	therefore,
	\[
		R(\phi) \geq R(\bar{\phi}) \geq \log_{2}(\log_{2}(n)).
	\]
	As this holds for all \(n \geq 1\) and \(\log_{2}(\log_{2}(n))\) tends to infinity as \(n\) increases, we must have that \(R(\phi) = \infty\).
\end{proof}
\begin{cor}	\label{cor:twistedStabilierFGTFResiduallyFinite}
	Let \(G\) be a torsion-free finitely generated residually finite group. Let \(\phi \in \Aut(G)\). If \(R(\phi) < \infty\), then all \(\phi\)-stabilizers are trivial.
\end{cor}
\begin{proof}
	As \(\phi\)-stabilizers are subgroups, a non-trivial one is necessarily infinite.
\end{proof}

The condition that \(G\) is finitely generated cannot be dropped. In fact, there already exists a non-finitely generated torsion-free residually finite abelian group admitting an automorphism with finite Reidemeister number and non-trivial twisted stabilizers.
\begin{example}[{Based on \cite[Proposition~3.6]{DekimpeGoncalves15}}]
	Consider the direct sum \(A\) of countably many copies of \(\Z\), indexed by the positive integers, \ie
	\[
		A = \bigoplus_{n = 1}^{\infty} \Z,
	\]
	and define \(\phi: A \to A\) as
	\[
		(a_{1}, a_{2}, a_{3}, a_{4}, \ldots) \mapsto (a_{1} + a_{2} + a_{3}, a_{2} + a_{3}, a_{3} + a_{4} + a_{5}, a_{4} + a_{5}, \ldots).
	\]
	In other words, the \((2k - 1)\)-th component is given by \(a_{2k - 1} + a_{2k} + a_{2k + 1}\) and the \(2k\)-th by \(a_{2k} + a_{2k + 1}\). This map is an endomorphism of \(A\), and even an automorphism since the map \(\psi: A \to A\) defined as
	\[
		(a_{1}, a_{2}, a_{3}, a_{4}, \ldots) \mapsto (a_{1} - a_{2}, a_{2} - a_{3} + a_{4}, a_{3} - a_{4}, a_{4} - a_{5} + a_{6}, \ldots)
	\]
	is the inverse of \(\phi\). The map \(\phi\) has non-trivial fixed points, for instance \((1, 0, 0, \ldots)\), but \(R(\phi) = 1\), since 
	\[
		[(0, 0, \ldots)]_{\phi} = \{a - \phi(a) \mid a \in A\} = \im(\phi - \Id)
	\]
	and \(\phi - \Id: A \to A\) is given by
	\[
		(a_{1}, a_{2}, a_{3}, a_{4}, \ldots) \mapsto (a_{2} + a_{3}, a_{3}, a_{4} + a_{5}, a_{5}, \ldots),
	\]
	which is clearly surjective
\end{example}

\begin{theorem}	\label{theo:SpecRDirectProductFGTFResiduallyFiniteCharacteristic}
	Let \(G\) and \(H\) be groups with \(G\) finitely generated, torsion-free residually finite. Suppose that \(H\) is characteristic in \(G \times H\). Then
	\[
		\SpecR(G \times H) = \SpecR(G) \cdot \SpecR(H).
	\]
\end{theorem}
\begin{proof}
	Since \(H\) is characteristic in \(G \times H\), each automorphism of \(G \times H\) is of the form
	\[
		\begin{pmatrix}
			\alpha & 0	\\ \gamma & \delta
		\end{pmatrix},
	\]
	and \(\alpha \in \Aut(G)\), \(\delta \in \Aut(H)\) and \(\gamma \in \Hom(G, Z(H))\), by \cref{lem:upperTriangularAutomorphisms}.
	Now, fix \(\phi \in \Aut(G \times H)\). We claim that \(R(\phi) = R(\alpha) R(\delta)\). We know by \cref{lem:ReidemeisterNumbersExactSequence}\eqref{lemitem:ReidemeisterNumberGreaterThanQuotientReidemeisterNumber} that \(R(\phi) \geq R(\alpha)\). Thus, \(R(\phi) = \infty\) if  \(R(\alpha) = \infty\). In that case, \(R(\phi) = R(\alpha)R(\delta)\). So suppose that \(R(\alpha) < \infty\). Then, by \cref{cor:twistedStabilierFGTFResiduallyFinite}, we know that \(\Stab_{\alpha}(g) = 1\) for all \(g \in G\). This means that the action of \(\Stab_{\alpha}(g)\) on \(\Reid[\delta]\) is trivial for all \(g \in G\), implying that the number of orbits for each action is equal to \(R(\delta)\). Applying \cref{theo:sumFormulaReidemeisterNumberDirectProductInvariantFactor} then yields
	\[
		R(\phi) = \sum_{[g]_{\alpha} \in \Reid[\alpha]} R(\delta) = R(\alpha) R(\delta).
	\]
	This proves that \(\SpecR(G \times H) \subseteq \SpecR(G) \cdot \SpecR(H)\). The converse inclusion follows directly from \cref{prop:ReidemeisterNumberDirectProductAutomorphismGroups}.
\end{proof}

\begin{cor}	\label{cor:SpecRDirectProductFGTFResiduallyFiniteFinite}
	Let \(G\) be a finitely generated, torsion-free residually finite group and \(F\) a group which is generated by torsion elements. Then
	\[
		\SpecR(G \times F) = \SpecR(G) \cdot \SpecR(F).
	\]
\end{cor}
\begin{proof}
	Since \(G\) is torsion-free, \(F\) is characteristic in \(G \times F\) as it is the subgroup generated by the torsion elements. Therefore, we can apply \cref{theo:SpecRDirectProductFGTFResiduallyFiniteCharacteristic}.
\end{proof}
Using the fact that finitely generated linear groups are residually finite (see \eg \cite{Malcev65}) and virtually torsion-free if the underlying field has characteristic zero (see \eg \cite{Selberg62}), one can easily generate examples of groups to which the theorem and its corollary apply. In particular, they apply to finitely generated torsion-free nilpotent groups and for several families of these groups their Reidemeister spectrum is known, see \eg \cite{DekimpeGoncalves14,DekimpeTertooyVargas20,Tertooy19}.

\section{Direct Products of Centerless Groups}
In this section we use the matrix description of the endomorphism monoid to describe the automorphism group under certain conditions. The first result is a generalization of F.\ Johnson's result \cite[Corollary~2.2]{Johnson83}.
\begin{defin}
	Let \(G\) be a group. We say that \(G\) is \emph{directly indecomposable} if \(G \cong H \times K\) for some groups \(H, K\) implies that \(H = 1\) or \(K = 1\).
\end{defin}
\begin{theorem}	\label{theo:generalisationJohnson83}
	Let \(G_{1}, \ldots, G_{n}\) be non-isomorphic non-trivial, centerless, directly indecomposable groups. Let \(r_{1}, \ldots, r_{n}\) be positive integers. Put \(G = \Times\limits_{i = 1}^{n} G_{i}^{r_{i}}\). Then
	\[
		\Aut(G) = \Times_{i = 1}^{n} \left(\Aut(G_{i}) \wr S_{r_{i}}\right)
	\]
\end{theorem}
In \cite{Johnson83}, there is the extra condition that there do not exist non-trivial homomorphisms \(\phi: G_{i} \to G_{j}\) with normal image for \(1 \leq i < j \leq n\), but it is redundant. In fact, the proof we will give is nearly identical to the one F.\ Johnson gave for the case \(n = 1\), see \cite[Theorem~1.1]{Johnson83}.

The following is well-known, so we omit the proof.
\begin{lemma}	\label{lem:equivalenceDirectProductNormalSubgroups}
	Let \(G\) be a group and let \(N_{1}, \ldots, N_{n}\) be commuting normal subgroups such that \(N_{1} \ldots N_{n} = G\). If
	\[
		\left(\prod_{j \ne i} N_{j}\right) \cap N_{i} = 1
	\]
	for all \(i \in \{1, \ldots, n\}\), then \(G \cong N_{1} \times \ldots \times N_{n}\). 
\end{lemma}

To make the notation easier, we first prove the following (equivalent) result:
\begin{theorem}	\label{theo:equivalentFormulationGeneralisationJohnson}
	Let \(G_{1}, \ldots, G_{n}\) be non-trivial, centerless, directly indecomposable groups. Put \(G = \Times\limits_{i = 1}^{n} G_{i}\). Under the monoid isomorphism \(\End(G) \cong \M\), \(\Aut(G)\) corresponds to those matrices \((\phi_{ij})_{ij}\) satisfying the following conditions:
	\begin{enumerate}[(i)]
		\item Each row and column contains exactly one non-trivial homomorphism.
		\item Each non-trivial homomorphism is an isomorphism.
	\end{enumerate}
\end{theorem}
\begin{proof}
	Let \(\phi = (\phi_{ij})_{ij} \in \Aut(G)\) and fix \(i \in \{1, \ldots, n\}\). By \cref{lem:automorphismOfDirectProductImpliesNormalImages}, we know that \(G_{i}\) is generated by \(\im \phi_{i1}\) up to \(\im \phi_{in}\) and that each of these images is normal in \(G_{i}\). Let \(k \in \{1, \ldots, n\}\) be arbitrary. Suppose that \(g\) is an element of
	\[
		\left(\prod_{j \ne k} \im \phi_{ij}\right) \cap \im \phi_{ik} \leq G_{i}.
	\]
For \(l \ne k\), we find that
	\[
		[g, \im \phi_{il}] \subseteq [\im \phi_{ik}, \im \phi_{il}] = 1,
	\]
	and for \(l = k\), we find that
	\[
		[g, \im \phi_{ik}] \subseteq \left[\prod_{j \ne k} \im \phi_{ij}, \im \phi_{ik}\right] = 1.
	\]
	Therefore, \(g \in Z(G_{i}) = 1\). Thus, we conclude that \(G_{i}\) is isomorphic to the direct product of \(\im \phi_{i1}\) up to \(\im \phi_{in}\) by \cref{lem:equivalenceDirectProductNormalSubgroups}. As \(G_{i}\) is directly indecomposable, exactly one of these images is non-trivial. Since \(i\) was arbitrary, we thus find a map \(\sigma: \{1, \ldots, n\} \to \{1, \ldots, n\}\) such that \(\im \phi_{i j} \ne 1\) if and only if \(j = \sigma(i)\).
	
	If \(\sigma\) were not surjective, say, \(m \notin \im \sigma\), we find that \(\im \phi_{im} = 1\) for all \(i \in \{1, \ldots, n\}\). Recall that \(e_{i}: G_{i} \to G\) and \(\pi_{i}: G \to G_{i}\) were the canonical inclusion and projection, respectively. Then \(e_{m}(G_{m}) \leq \ker (\pi_{i} \circ \phi)\) for all \(i \in \{1, \ldots, n\}\). We thus find that
	\[
		e_{m}(G_{m}) \leq \bigcap_{i = 1}^{n} \ker(\pi_{i} \circ \phi) = \bigcap_{i = 1}^{n} \inv{\phi}(\ker(\pi_{i})) =  \inv{\phi} \left(\bigcap_{i = 1}^{n} \ker (\pi_{i}) \right) = 1,
	\]
	which contradicts the non-triviality of \(G_{m}\). Thus, \(\sigma\) is surjective, and therefore bijective. We thus find that the matrix representation of \(\phi\) contains exactly one non-trivial homomorphism on each row and column. As \(\phi\) is an automorphism, each of these non-trivial homomorphisms must be both injective and surjective, therefore, they must all be isomorphisms.
\end{proof}
\begin{proof}[Proof of \cref{theo:generalisationJohnson83}]
	Let \(\phi \in \Aut(G)\). By \cref{theo:equivalentFormulationGeneralisationJohnson}, the matrix representation \((\phi_{ij})_{ij}\) contains exactly one non-trivial homomorphism per row and column, and each of those homomorphisms is in fact an isomorphism. Due to the fact that \(G_{i}\) and \(G_{j}\) are not isomorphic for \(i \ne j\), \((\phi_{ij})_{ij}\) is of the form
	\[
		\begin{pmatrix}
			A_{1}	&	0		&	0		&	\ldots	&	0	\\
			0		&	A_{2}	&	0		&	\ldots	&	0	\\
			0		&	0		&	A_{3}	&	\ldots	&	0	\\
			\vdots	&	\vdots	&	\vdots	&	\ddots	&	\vdots	\\
			0		&	0		&	0		&	\ldots	&	A_{n},
		\end{pmatrix}
	\]
	where each \(A_{i}\) is an \((r_{i} \times r_{i})\)-matrix containing exactly one non-trivial homomorphism per row and column, and each of those homomorphisms is an automorphism of \(G_{i}\). These block matrices correspond to automorphisms lying in \(\Aut(G_{i}) \wr S_{r_{i}}\). Therefore, \(\phi\) lies in \(\Times\limits_{i = 1}^{n} \left(\Aut(G_{i}) \wr S_{r_{i}}\right)\).
	\end{proof}

\begin{cor}	\label{cor:SpecRProductCentrelessDirectlyIndecomposableGroups}
	Let \(G_{1}, \ldots, G_{n}\) be non-trivial, non-isomorphic, centerless, directly indecomposable groups. Let \(r_{1}, \ldots, r_{n}\) be positive integers and put \(G = \Times\limits_{i = 1}^{n}G_{i}^{r_{i}}\). Then
	\[
		\SpecR(G) = \prod_{i = 1}^{n} \left(\bigcup_{j = 1}^{r_{i}} \SpecR(G_{i})^{(j)}\right).
	\]
	In particular, \(G\) has the \(\Rinf\)-property if and only if \(G_{i}\) has the \(\Rinf\)-property for some \(i \in \{1, \ldots, n\}\).
\end{cor}
\begin{proof}
	This follows by combining \cref{theo:generalisationJohnson83} with \cref{cor:SpecRDirectProductCharacteristicSubgroups,cor:ReidemeisterNumbersOfAut(G)wrSn}.
\end{proof}
\begin{remark}
	S.\ Tertooy proved the same equality for direct products of free nilpotent groups of finite rank in \cite[\S6.5]{Tertooy19}. Note that free nilpotent groups have non-trivial center, so they do not satisfy the conditions of \cref{cor:SpecRProductCentrelessDirectlyIndecomposableGroups}, showing that those conditions are sufficient but not necessary.
\end{remark}

\begin{example}	\label{ex:directProductsOfFreeProducts}
	Recall the result by D.\ Gon\c{c}alves, P.\ Sankaran and P.\ Wong \cite[Theorem~1]{GoncalvesSankaranWong20} stating that, under some mild conditions on the \(G_{i}\), the free product \(G_{1} * \ldots * G_{n}\) with \(n \geq 2\) has the \(\Rinf\)-property. Since free products are centerless (see  \eg \cite[Corollary~4.5]{KarrassMagnusSolitar76}) and directly indecomposable (see \eg \cite[Observation~p.~177]{LyndonSchupp77}), \cref{cor:SpecRProductCentrelessDirectlyIndecomposableGroups} then implies that any direct product of such free products has the \(\Rinf\)-property. This includes, for instance, direct products of groups isomorphic to free groups of finite rank, the modular group \(\PSL(2, \Z) \cong \Z/2\Z * \Z / 3\Z\) or the infinite dihedral group \(D_{\infty} \cong \Z / 2 \Z * \Z / 2 \Z\).
\end{example}

\begin{example}
	Several known results regarding automorphism groups and the \(\Rinf\)-property follow from \cref{theo:generalisationJohnson83} and \cref{cor:SpecRProductCentrelessDirectlyIndecomposableGroups}. For instance, A.\ Fel'shtyn and T.\ Nasybullov proved in \cite[Theorem~3]{FelshtynNasybullov16} that certain reductive linear algebraic groups \(G\) have the \(\Rinf\)-property by proving it for the quotient group \(G / R(G)\), which splits as a direct product of Chevalley groups. The latter can also be proved by combining the results from \cite{Nasybullov20,Nasybullov16,Nasybullov15a,Nasybullov12,Nasybullov19} with \cref{cor:SpecRProductCentrelessDirectlyIndecomposableGroups}.
	
	Another example concerns right-angled Artin groups. Each RAAG admits a unique maximal decomposition as a direct product of RAAGs (see \eg \cite[Proposition~3.1]{GandiniWahl16}) and N.\ Fullarton \cite{Fullarton16} and G.\ Gandini and N.\ Wahl \cite{GandiniWahl16} described the automorphism group of the RAAG in terms of the automorphism groups of the direct factors. For centerless RAAGs, this description also follows from \cref{theo:generalisationJohnson83}.
%
\end{example}

\begin{example}	\label{ex:directProductsCrystallographicGroups}
	For each \(n \geq 2\), let \(H_{n} := \Z^{n} \rtimes_{-I} \Z/2\Z\), where the action is given by inversion. K.\ Dekimpe, T.\ Kaiser and S.\ Tertooy determined the Reidemeister spectra of these groups \cite[Proposition~5.7]{DekimpeKaiserTertooy19}:
	\[
		\SpecR(H_{n}) = \begin{cases}
			2\N_{0} \cup \{3, \infty\}	&	\mbox{if } n = 2	\\
			\N_{0} \setminus \{1\} \cup \{\infty\}	&	\mbox{if } n \geq 3.
		\end{cases}
	\]
	Moreover, we argue that each \(H_{n}\) is directly indecomposable and centerless as well. We write the elements of \(H_{n}\) as \((x, t^{i})\), where \(x \in \Z^{n}\) and \(i \in \{0, 1\}\). Consider an element of the form \((x, t)\). For \(y \in \Z^{n}\) and \( i \in \{0, 1\}\), note that
	\begin{align*}
		(y, t^{i})(x, t) \inv{(y, t^{i})} 	&= (y + (-1)^{i}x, t^{i + 1}) ((-1)^{i + 1}y, t^{i})	\\
							&= (y + (-1)^{i}x + y, t) = ((-1)^{i}x + 2y, t).
	\end{align*}
	For this to be equal to \((x, t)\), we must have \(y = 0\) and \(i = 0\), or \(y = x\) and \(i = 1\). Consequently, the centralizer \(C_{H_{n}}(x, t)\) of \((x, t)\) equals \(\grpgen{(x, t)} \cong \Z/ 2\Z\). If we take two distinct \(x, y \in \Z^{n}\), then \(C_{H_{n}}(x, t) \cap C_{H_{n}}(y, t) = 1\). Since the center of \(H_{n}\) is contained in this intersection, the center is trivial.
	
	Now, suppose that \(H_{n} \cong A \times B\). Under the isomorphism, one of the factors, say \(A\), must contain the image \(g\) of some \((x, t)\). The centralizer of \(g\) is then isomorphic to \(\Z/2\Z\), but it also equals \(C_{A}(g) \times B\). Since \(C_{A}(g)\) is non-trivial and \(\Z / 2\Z\) is directly indecomposable, it must hold that \(B = 1\). Therefore, \(H_{n}\) is directly indecomposable.
	
	Consequently, we can apply \cref{cor:SpecRProductCentrelessDirectlyIndecomposableGroups} to compute \(\SpecR(H_{n_{1}}^{r_{1}} \times \ldots \times H_{n_{k}}^{r_{k}})\) for all \(2 \leq n_{1} < \ldots < n_{k}\) and \(r_{1}, \ldots, r_{k} \geq 1\). For instance,
	\[
		\SpecR(H_{2} \times H_{3}) = \big((2 \N_{0} \cup 3 \N_{0}) \setminus \{2, 3\} \big) \cup \{\infty\}.
	\]
\end{example}

The next result provides us with sufficient conditions for a direct product to have the \(\Rinf\)-property when one of the factors has it.
\begin{theorem}	\label{theo:SecondGeneralisationJohnson}
	Let \(G_{1}, \ldots, G_{n}, H\) be non-trivial groups such that each \(G_{i}\) is centerless and directly indecomposable and such that, for each \(i \in \{1, \ldots, n\}\), \(H\) has no direct factor isomorphic with \(G_{i}\). Let \(G = \Times\limits_{i = 1}^{n} G_{i}\). Then, under the monoid isomorphism \(\End(G \times H) \cong \M\), \(\Aut(G \times H)\) corresponds to
	\[
		\left\{ \begin{pmatrix} \alpha & 0 \\ \gamma & \delta \end{pmatrix} \middlebar \begin{array}{ccc}\alpha \in \Aut(G), &\delta \in \Aut(H), \\ \gamma \in \Hom(G, Z(H)) \end{array}\right\}.
	\]
	In other words, \(H\) is characteristic in \(G \times H\).
\end{theorem}
\begin{proof}
Let \(\phi \in \Aut(G \times H)\) and write
\[
	\phi = \begin{pmatrix}
		\phi_{11}		&	\ldots	&	\phi_{1n}	&	 \beta_{1}	\\
		\vdots		&	\ddots	&	\vdots		&	\vdots	\\
		\phi_{n1}		&	\ldots	&	\phi_{nn}		&	\beta_{n}	\\
		\gamma_{1}	&	\ldots	&	\gamma_{n}		&	\delta
	\end{pmatrix}
\]
with, for all \(1 \leq i, j \leq n\), \(\phi_{ij} \in \Hom(G_{j}, G_{i}), \beta_{i} \in \Hom(H, G_{i}), \gamma_{j} \in \Hom(G_{j}, H)\) and \(\delta \in \End(H)\). By a similar argument as in the beginning of \cref{theo:equivalentFormulationGeneralisationJohnson}, we find that each of the first \(n\) rows of \(\phi\) contains precisely one non-trivial homomorphism. The inverse automorphism \(\inv{\phi}\) has a similar matrix form as \(\phi\), hence, also each of the first \(n\) rows of \(\inv{\phi}\) contains precisely one non-trivial homomorphism, and we denote the corresponding homomorphisms of \(\inv{\phi}\) by \(\psi_{ij}, \beta_{i}', \gamma_{j}'\) and \(\delta'\).

The goal is to prove that \(\beta_{i} = 0\) for all \(1 \leq i \leq n\). First, we make some observations.

Suppose that \(\phi_{ij} \ne 0\) for some \(i, j \in \{1, \ldots, n\}\). Then \(\beta_{i} = 0\) and \(\phi_{ik} = 0\) for \(k \in \{1, \ldots, n\}\) different from \(j\). Hence \(\Id_{G_{i}} = (\phi \circ \inv{\phi})_{ii} = \phi_{ij} \psi_{ji}\). This implies that \(\psi_{ji} \ne 0\), hence by symmetry, \(\Id_{G_{j}} = (\inv{\phi} \circ \phi)_{jj} = \psi_{ji} \phi_{ij}\). Moreover, this also shows that \(\phi_{kj} = 0\) for \(k \ne i\). Indeed, if \(\phi_{kj} \ne 0\), the same arguments as before show that \(\psi_{jk} \ne 0\). This would imply that the \(j\)-th row of \(\inv{\phi}\) contains two non-trivial homomorphisms, which is a contradiction. Hence, the index \(i\) is unique, \ie the \(j\)-th column of \(\phi\) contains only \(\phi_{ij}\) and \(\gamma_{j}\) as (potentially) non-trivial homomorphisms. By symmetry, the \(i\)-th column of \(\inv{\phi}\) contains only \(\gamma_{i}'\) and \(\psi_{ji}\) as (potentially) non-trivial homomorphisms.

Next, suppose that \(\beta_{l}\) is non-trivial for some \(l \in \{1, \ldots, n\}\). Since each of the first \(n\) rows of \(\phi\) contains at most one non-trivial homomorphism, there is a column of \(\phi\), say the \(j\)-th one, containing only \(\gamma_{j}\) as (potentially) non-trivial homomorphism. So, let \(\J\) be the set of indices of the columns of \(\phi\) of this form and \(\I\) be the set of indices \(i\) with \(\beta_{i} \ne 0\). Suppose that \(\I\), and thus also \(\J\), is non-empty.

Note that, for each \(j \in \J\), \(\gamma_{j}\) is injective, since
\[
	\Id_{G_{j}} = (\inv{\phi} \circ \phi)_{jj} = \beta_{j} \gamma_{j}.
\]
From this it follows that \(\im \gamma_{j} \cap \grpgen{\im \gamma_{i} \mid i \in \J \setminus \{j\}} = 1\) for \(j \in \J\). Indeed, suppose that \(h\) lies in the intersection. Since the images of the \(\gamma_{k}\)'s commute pairwise, we can write
\[
	h = \prod_{i \in \J \setminus \{j\}} \gamma_{i}(g_{i}) = \gamma_{j}(g_{j})
\]
for some \(g_{j} \in G_{j}\) and \(g_{i} \in G_{i}\). For \(g \in G_{j}\), we then find that
\[
	\gamma_{j}(g) \gamma_{j}(g_{j}) = \gamma_{j}(g) \prod_{i \in \J \setminus \{j\}}\gamma_{i}(g_{i}) = \left(\prod_{i \in \J \setminus \{j\}}\gamma_{i}(g_{i})\right) \gamma_{j}(g) = \gamma_{j}(g_{j}) \gamma_{j}(g),
\]
again by the commuting condition. Injectivity of \(\gamma_{j}\) implies that \(gg_{j} = g_{j}g\), and as this has to hold for all \(g \in G_{j}\), we conclude that \(g_{j} \in Z(G_{j}) = 1\). Thus, \(h = 1\), proving the claim. This shows that \(\grpgen{\im \gamma_{j} \mid j \in \J}\) is isomorphic to \(\Times\limits_{j \in \J} G_{j}\), by \cref{lem:equivalenceDirectProductNormalSubgroups}. We thus have a subgroup of \(H\) isomorphic to \(\Times\limits_{j \in \J} G_{j}\). We now proceed to prove that this subgroup is in fact a direct factor to obtain a contradiction.

First, we claim that \(\delta \gamma'_{i} = 0\) for \(i \in \I\). In order to do so, we compute \({(\phi \circ \inv{\phi})_{n+1, i}}\) and obtain
\[
	0 = \delta \gamma'_{i} + \sum_{j = 1}^{n} \gamma_{j} \psi_{ji}.
\]
If \(\psi_{ji}\) would be non-trivial, then also \(\phi_{ij}\) would be non-trivial, which would yield two non-trivial homomorphisms on the \(i\)-th row of \(\phi\). Hence, \(0 = \delta \gamma_{i}'\), proving the claim.

With this equality, we can prove that \(\delta \delta' \delta = \delta\), by noting that
\[
	\Id_{H} = (\inv{\phi} \circ \phi)_{n + 1, n + 1} = \delta' \delta + \sum_{i \in \I} \gamma_{i}' \beta_{i},
\]
as \(\I\) contains all indices \(i\) for which \(\beta_{i}\) is non-trivial. Composing with \(\delta\) and using \(\delta \gamma_{i}' = 0\) yields \(\delta = \delta \delta' \delta\) as desired.

Next, we prove the following equality for all \(j \in \{1, \ldots, n\}\):
\begin{equation}	\label{eq:deltadelta'gammam}
	\delta \delta' \gamma_{j} = \begin{cases}
		0	&	\mbox{if } j \in \J	\\
		\gamma_{j}	&	\mbox{if } j \notin \J
	\end{cases}
\end{equation}
Fix \(j \in \{1, \ldots, n\}\). Suppose that \(j \in \J\). Then \(\phi_{ij} = 0\) for all \(i \in \{1, \ldots, n\}\). Consequently, \(0 = (\inv{\phi} \circ \phi)_{n + 1, j} = \delta' \gamma_{j}\), hence \(\delta \delta' \gamma_{j} = 0\) as well.

If \(j \notin \J\), then there is a unique \(i \in \{1, \ldots, n\}\) such that \(\phi_{ij} \ne 0\). Then also \(\psi_{ji} \ne 0\), and this will be the only non-trivial homomorphism in the first \(n\) rows of the \(i\)-th column of \(\inv{\phi}\). Using this, we find that
\[
	0 = (\inv{\phi} \circ \phi)_{n + 1, j} = \gamma_{i}' \phi_{ij} + \delta' \gamma_{j}
\]
and
\[
	0 = (\phi \circ \inv{\phi})_{n + 1, i} = \gamma_{j} \psi_{ji} + \delta \gamma_{i}'.
\]
Composing the first equality on the left with \(\delta\) yields
\begin{align*}
	0	&= \delta \gamma_{i}' \phi_{ij} + \delta \delta' \gamma_{j}	\\
		&= - \gamma_{j} \psi_{ji} \phi_{ij} + \delta \delta' \gamma_{j}	\\
		&= - \gamma_{j} + \delta \delta' \gamma_{j},
\end{align*}
where we used that \(\Id_{G_{j}} = \psi_{ji} \phi_{ij}\). Rearranging terms proves that \(\delta \delta' \gamma_{j} = \gamma_{j}\), finishing the proof of \eqref{eq:deltadelta'gammam}.

Finally, we prove that
\begin{equation}	\label{eq:firstTrivialIntersectionDirectProduct}
	\grpgen{\{\im \gamma_{j} \mid j \in \J\}} \cap \grpgen{\im \delta, \{\im \gamma_{j} \mid j \notin \J\}} = 1.
\end{equation}
Let \(h\) be an element of the intersection. Since all the images of the \(\gamma_{i}\)'s and \(\delta\) commute pairwise, we can write
\[
	h = \prod_{j \in \J} \gamma_{j}(g_{j}) = \delta(h') \prod_{j \notin \J} \gamma_{j}(g_{j})
\]
where \(h' \in H\) and \(g_{j} \in G_{j}\) for all indices \(j\). Applying \(\delta \delta'\) on this equality and using \eqref{eq:deltadelta'gammam} yields on the one hand
\[
	\delta \delta'(h) = \prod_{j \in \J} \delta \delta' \gamma_{j}(g_{j}) = 1,
\]
and on the other hand, again by using \eqref{eq:deltadelta'gammam} and \(\delta = \delta \delta' \delta\),
\[
	\delta \delta'(h) = \delta \delta' \delta (h') \prod_{j \notin \J} \delta \delta' \gamma_{j}(g_{j}) = \delta(h')\prod_{j \notin \J} \gamma_{j}(g_{j}) = h.
\]
Hence, \(h = 1\).

Now, both subgroups in \eqref{eq:firstTrivialIntersectionDirectProduct} are normal, being a product of normal subgroups by \cref{lem:automorphismOfDirectProductImpliesNormalImages}, they generate \(H\), also by \cref{lem:automorphismOfDirectProductImpliesNormalImages} and they intersect trivially. Therefore, we conclude that \(H\) equals the internal direct product of the groups \(\im \phi_{j}\) with \(j \in \J\) and the group \(\grpgen{\im \delta, \{\im \gamma_{j} \mid j \notin \J\}}\), \ie
\[
	H \cong \left(\Times_{j \in \J} G_{j}\right) \times N
\]
for some normal subgroup \(N\) of \(H\). As \(|\J| \geq 1\), we obtain a contradiction, since we assume that \(H\) has no direct factors isomorphic with \(G_{j}\) for \(j \in \{1, \ldots, n\}\). Therefore, \(\J\) is empty, hence also \(\I\) is empty, meaning that all \(\beta_{i}\)'s are trivial. Thus, the matrix form of \(\phi\) is
\[
	\phi = \begin{pmatrix}
		\phi_{11}		&	\ldots	&	\phi_{1n}	&	 0	\\
		\vdots		&	\ddots	&	\vdots		&	\vdots	\\
		\phi_{n1}		&	\ldots	&	\phi_{nn}		&	0	\\
		\gamma_{1}	&	\ldots	&	\gamma_{n}		&	\delta
	\end{pmatrix}.
\]
Rewrite this as
\[
	\phi = \begin{pmatrix}
		\alpha	&	0	\\
		\gamma	&	\delta
	\end{pmatrix}
\]
with \(\alpha \in \End(G), \gamma \in \Hom(G, H), \delta \in \End(H)\). As \(\phi\) was arbitrary, each automorphism of \(G \times H\) is of this form. \cref{lem:upperTriangularAutomorphisms} then implies that \(\alpha \in \Aut(G)\), \(\delta \in \Aut(H)\) and \(\gamma \in \Hom(G, Z(H))\).
\end{proof}

\begin{remark}
	By \cref{theo:equivalentFormulationGeneralisationJohnson}, we also know what \(\Aut(G)\) looks like. In fact, we could have merged \cref{theo:equivalentFormulationGeneralisationJohnson,theo:SecondGeneralisationJohnson} into one theorem, since the start of the proof of \cref{theo:SecondGeneralisationJohnson} is the same as that of \cref{theo:equivalentFormulationGeneralisationJohnson}. However, for the sake of clarity, we have split the results: \cref{theo:equivalentFormulationGeneralisationJohnson} deals with the internal structure of \(\Aut(G)\) and \cref{theo:SecondGeneralisationJohnson} determines the influence of \(G\) on an additional factor \(H\) on which we impose less strict conditions.
\end{remark}

\begin{cor}	\label{cor:SpecRDirectProductsH}
	Let \(G_{1}, \ldots, G_{n}, H\) and \(G\) be as in the previous theorem.
	\begin{enumerate}[(i)]
		\item If \(G_{i} \in \Rinf\) for some \(i \in \{1, \ldots, n\}\), then \(G \times H \in \Rinf\) as well.
		\item If \(G_{i}\) is finitely generated residually finite for each \(i \in \{1, \ldots, n\}\) and \(H \in \Rinf\), then \(G \times H \in \Rinf\) as well.
		\item If \(G_{i}\) is finitely generated torsion-free residually finite for each \(i \in \{1, \ldots, n\}\), then
		\[
			\SpecR(G \times H) = \SpecR(G) \cdot \SpecR(H).
		\]
	\end{enumerate}
\end{cor}
\begin{proof}
	By the previous theorem, \(H\) is characteristic in \(G \times H\). If \(G_{i} \in \Rinf\) for some \(i \in \{1, \ldots, n\}\), then so does \(G\) by \cref{cor:SpecRProductCentrelessDirectlyIndecomposableGroups} and therefore \(G \times H\) has the \(\Rinf\)-property as well by \cref{lem:ReidemeisterNumbersExactSequence}, since \(G \cong \frac{G \times H}{H}\).
	
	Now, suppose that \(H \in \Rinf\) and that each \(G_{i}\) is finitely generated residually finite. Then also \(G\) is finitely generated residually finite. For an automorphism \(\phi\) of \(G \times H\), denote by \(\phi'\) the induced automorphism on \(H\) and by \(\bar{\phi}\) the induced automorphism on \(G\). If \(R(\bar{\phi}) = \infty\), then \cref{lem:ReidemeisterNumbersExactSequence}\eqref{lemitem:ReidemeisterNumberGreaterThanQuotientReidemeisterNumber} yields \(R(\phi) = \infty\). If \(R(\bar{\phi}) < \infty\), then \cref{prop:infiniteTwistedStabilisersImpliesInfiniteReidemeisterNumberFGRF} implies that \(\Fix(\bar{\phi})\) is finite. Since \(H \in \Rinf\), we know that \(R(\phi') = \infty\) and then \cref{lem:ReidemeisterNumbersExactSequence}\eqref{lemitem:InfiniteSubgroupReidemeisterNumberFiniteFixedPointsImpliesInfiniteReidemeisterNumber} yields that also \(R(\phi) = \infty\). 
	
	The last item follows from \cref{theo:SpecRDirectProductFGTFResiduallyFiniteCharacteristic}.
\end{proof}

\begin{example}
We continue example \cref{ex:directProductsOfFreeProducts}. Let \(G_{1}\) up to \(G_{n}\) each be a non-trivial free product of groups satisfying the same conditions as in \cref{ex:directProductsOfFreeProducts} and let \(H\) be a finite group. Then each \(G_{i}\) has the \(\Rinf\)-property and since a (non-trivial) free product is infinite, it cannot be a direct factor of a finite group. Thus, the direct product \(G_{1} \times \ldots \times G_{n} \times H\) has the \(\Rinf\)-property as well by the first item of \cref{cor:SpecRDirectProductsH}.

\end{example}

\begin{example}
We continue \cref{ex:directProductsCrystallographicGroups}. For each \(k \geq 2\), the group \(H_{k}\) contains torsion and since \(H_{k}\) fits in the exact sequence \(1 \to \Z^{k} \to H_{k} \to \Z/2 \Z \to 1\), it is a finitely generated residually finite group. On the other hand, the non-abelian Baumslag-Solitar groups \(BS(m, n)\) have the \(\Rinf\)-property \cite[Theorem~4.4]{FelshtynGoncalves06}, and are torsion-free, so they cannot contain \(H_{k}\) as a direct factor. Therefore, if we consider direct products of the form \(\left(\Times\limits_{i = 1}^{k} H_{n_{i}}^{r_{i}}\right) \times BS(m, n)\) for some \(2 \leq n_{1} < \ldots < n_{k}\), \(r_{1}, \ldots, r_{k} \geq 1\) and \(BS(m, n)\) non-abelian, we can apply the second item of \cref{cor:SpecRDirectProductsH} to prove that these groups have the \(\Rinf\)-property.
\end{example}

\begin{example}
This example is based on the results from \cite{GoncalvesWong03}. Consider the semi-direct product \(G = \Z^{2} \rtimes_{A} \Z\), where \(A = \begin{pmatrix} 2 & 1 \\ 1 & 1\end{pmatrix}\). As in \cref{ex:directProductsCrystallographicGroups}, we write the elements in the group as \((x, t^{i})\), where \(x \in \Z^{2}\) and \(i \in \Z\). It is a polycyclic group of Hirsch length \(3\), hence finitely generated and residually finite, and clearly torsion-free. Also note that \(G\) is solvable.

To prove that its center is trivial, let \((x, t^{i})\) be an element in the center. Conjugating with \((0, t)\) yields the equality \((x, t^{i}) = (Ax, t^{i})\). Since the eigenvalues of \(A\) are \(\frac{3 \pm \sqrt{5}}{2}\), \(A\) does not have any non-trivial fixed points, therefore, \(x = 0\). If we then conjugate \((0, t^{i})\) with \((y, 1)\), where \(y \in \Z^{2}\) is a non-trivial element, we find that \((0, t^{i}) = (y - A^{i}y, t^{i})\). Since neither of the eigenvalues of \(A\) is a root of unity, \(A^{i}\) has no eigenvalue \(1\) if \(i \ne 0\). Therefore, \(i = 0\) and we find that \((x, t^{i}) = (0, 1)\), proving that the center of \(G\) is trivial.

To prove that \(G\) is directly indecomposable, note that, for \(x \in \Z^{2}\), we have \([(x, 1), (0, \inv{t})] = (Ax - x, 1)\). Since \(\det(A - I) = -1\), \(A - I\) defines a surjective map on \(\Z^{2}\), thus \([G, G]\) contains \(\Z^{2}\). The inclusion \([G, G] \subseteq \Z^{2}\) follows immediately from the definition of \(G\). So, \([G, G] = \Z^{2}\) and \(G / [G, G] \cong \Z\). Next, suppose that \(G \cong H \times K\). Factoring out the commutator subgroup, we get the isomorphism \(\Z \cong \frac{H}{[H, H]} \times \frac{K}{[K, K]}\). As \(\Z\) is directly indecomposable, one of the factors, say \(H / [H, H]\), is trivial. This implies \([H, H] = H\). But \(H\) is isomorphic to a subgroup of \(G\), which is a solvable group, hence \(H\) is solvable itself. Combined with the equality \(H = [H, H]\), this yields \(H = 1\), showing that \(G\) is directly indecomposable.

Finally, in \cite[Theorem~4.1, Example~4.3]{GoncalvesWong03}, it is proven that \(G\) admits automorphisms with finite Reidemeister numbers. A concrete example of an automorphism with finite Reidemeister number is
\[
	\phi: G \to G: (x, t^{i}) \mapsto (Mx, t^{-i}),
\]
where \(M = \begin{pmatrix} 0 & 1 \\ -1 & 0 \end{pmatrix}\). One can verify that \(R(\phi) = 4\). Thus, if we consider the direct product \(G \times P\) where \(P\) is any polycyclic group of Hirsch length at most \(2\), we can apply the third item of \cref{cor:SpecRDirectProductsH} to conclude that \(\SpecR(G \times P) = \SpecR(G) \cdot \SpecR(P)\), since \(G\) cannot be a direct factor (not even a subgroup) of \(P\).
\end{example}


\section{Direct Products of Virtually Free Groups}
We now use the obtained results to investigate the \(\Rinf\)-property for direct products of finitely generated virtually free groups.

\begin{defin}
	Let \(G\) be a group. We say that \(G\) is \emph{non-elementary} virtually free if there is a non-abelian free subgroup \(F\) of finite index in \(G\).
\end{defin}
\begin{prop}	\label{prop:FGNonelementaryVirtuallyFreeRinf}
	Let \(G\) be a finitely generated non-elementary virtually free group. Then \(G\) has the \(\Rinf\)-property.
\end{prop}
\begin{proof}
	Since \(G\) is finitely generated and contains a non-abelian free subgroup of finite index, it is Gromov hyperbolic. The result then follows from \cite[Theorem~3]{Felshtyn01}.
\end{proof}

\begin{lemma}
	Let \(G\) be a non-elementary virtually free group. Then \(Z(G)\) is finite.
\end{lemma}
\begin{proof}
	Let \(F\) be a finite index non-abelian free subgroup of \(G\). If \(Z(G)\) is infinite, then \(F \cap Z(G)\) is non-trivial. However, \(F \cap Z(G)\) lies in the center of \(F\), which is trivial, since \(F\) is non-abelian free. Therefore, \(Z(G)\) is finite.
\end{proof}

\begin{lemma}	\label{lem:maximalFiniteNormalSubgroupVirtuallyFree}
	Let \(G\) be a virtually free group. Then \(G\) has a unique maximal finite normal subgroup \(N_{0}\), which is characteristic.
	
	Moreover, \(G / N_{0}\) has no non-trivial finite normal subgroups. In particular, if \(G\) is non-elementary virtually free, then \(G / N_{0}\) is centerless.
\end{lemma}
\begin{proof}
	Let \(F\) be a finite index free subgroup of \(G\). We may assume its normal, by intersecting all conjugates of \(F\). Denote by \(\pi: G \to G / F\) the canonical projection. Let \(E\) be a finite subgroup. Then \(E \cap F\) is trivial, since \(E\) is torsion. Therefore, \(\size{E} = \size{\pi(E)} \leq [G : F]\). Thus, the size of finite subgroups of \(G\) is bounded by \([G : F]\). Let \(N_{0}\) be a finite normal subgroup of maximal size. Let \(N\) be a arbitrary finite normal subgroup of \(G\). Then \(NN_{0}\) is a normal finite subgroup of \(G\). Therefore, \(\size{NN_{0}} \leq \size{N_{0}}\). Also, \(N_{0} \leq N N_{0}\), thus \(N_{0} = N N_{0}\), implying that \(N \leq N_{0}\). Thus, \(N_{0}\) is the unique maximal finite normal subgroup.
	
	Since automorphisms preserve order and normality, \(N_{0}\) is characteristic.
	
	Now, if \(N / N_{0}\) is a finite normal subgroup of \(G / N_{0}\), then \(N\) is a finite normal subgroup of \(G\) containing \(N_{0}\). Therefore, \(N = N_{0}\).
	
	For the center of \(G / N_{0}\), note that the (non-abelian) free subgroup \(F\) of \(G\) projects injectively down to \(G / N_{0}\). Therefore, \(G / N_{0}\) is non-elementary virtually free, hence has finite center. Since the center is normal, we conclude that \(Z(G / N_{0}) = 1\).
\end{proof}

\begin{lemma}	\label{lem:FiniteCharacteristicSubgroupVirtuallyFree}
	Let \(G_{1}, \ldots, G_{n}\) be non-elementary virtually free groups. Let \(N_{1}, \ldots, N_{n}\) be their respective maximal normal finite subgroups. Then \(N := \Times\limits_{i = 1}^{n}N_{i}\) is characteristic in \(G := \Times\limits_{i = 1}^{n} G_{i}\).
\end{lemma}
\begin{proof}
	Let \(\phi \in \Aut(G)\) and \(g = (g_{1}, \ldots, g_{n}) \in \phi(N)\). Let \(\pi_{j}: G \to G_{j}\) be the projection. Then \(\pi_{j}(\phi(N))\) is a finite normal subgroup of \(G_{j}\), therefore \(g_{j} \in \pi_{j}(\phi(N)) \leq N_{j}\). This holds for all \(j \in \{1, \ldots, n\}\), implying that \(g \in N\).
\end{proof}

\begin{lemma}	\label{lem:VirtuallyFreeGroupsDirectlyIndecomposable}
	Let \(G\) be a virtually free group with no non-trivial finite normal subgroups. Then \(G\) is directly indecomposable.
\end{lemma}
\begin{proof}
	Suppose that \(G \cong H \times K\) for some non-trivial subgroups of \(G\). Since \(H\) and \(K\) must be normal in \(G\), both are infinite. Let \(F\) be a free subgroup of finite index in \(G\). Then \(H \cap F\) and \(K \cap F\) are finite index subgroups of \(H\) and \(K\) as well, hence they are infinite. Therefore, \(F\) contains the subgroup \((H \cap F) \times (K \cap F)\). This is a contradiction, however, since \((H \cap F) \times (K \cap F)\) must be a free group (being a subgroup of \(F\)), but free groups are directly indecomposable (see \eg \cite[Observation~p.~177]{LyndonSchupp77})
\end{proof}

\begin{theorem}
	Let \(G_{1}, \ldots, G_{n}\) be finitely generated non-elementary virtually free groups. Then
	\(
		G := \Times\limits_{i = 1}^{n} G_{i}
	\)
	has the \(\Rinf\)-property.
\end{theorem}
\begin{proof}
	Let \(N\) be the characteristic subgroup of \(G\) provided by \cref{lem:FiniteCharacteristicSubgroupVirtuallyFree}. By \cref{lem:ReidemeisterNumbersExactSequence}, it is sufficient to prove that \(G / N\) has the \(\Rinf\)-property. Therefore, by \cref{lem:maximalFiniteNormalSubgroupVirtuallyFree}, we may assume that each \(G_{i}\) is centerless and has no non-trivial finite normal subgroups. \cref{lem:VirtuallyFreeGroupsDirectlyIndecomposable} then implies that each \(G_{i}\) is directly indecomposable. Since each \(G_{i}\) has the \(\Rinf\)-property by \cref{prop:FGNonelementaryVirtuallyFreeRinf}, \cref{cor:SpecRProductCentrelessDirectlyIndecomposableGroups} implies that \(G\) has it as well.	
\end{proof}
\section*{Acknowledgements}
The author thanks Karel Dekimpe and the anonymous referee for their useful remarks and suggestions, in particular for their examples of concrete groups to which the results apply.
\section*{Funding}
This work was supported by the Research Foundation -- Flanders (FWO) under Grant 1112520N.
\printbibliography
\end{document}